\documentclass{article}
\usepackage[all]{xy}
\usepackage{latexsym,amssymb,amsfonts,amsmath}
\usepackage{graphicx}
  \newtheorem{definition}{Def.}[section]
  \newtheorem{theorem}[definition]{Theorem}
  \newtheorem{lemma}[definition]{Lemma}
  \newtheorem{cor}[definition]{Corollary}
  \newtheorem{prop}[definition]{Proposition}
  \newtheorem{remark}[definition]{Remark}

\def\releps{\mathrel\epsilon} 

\newcommand{\primconst}[1]              
{{\mathop{\sf #1}}} 
\newcommand{\ext}{{\primconst{\,ext\,}}}
      
\newcommand{\rest}{\primconst{\,rest\,}}

\newcommand{\bpr}{{\bf Proof.} }

\def\releps{\mathrel\epsilon}

 \CompileMatrices

\makeatletter
\providecommand*{\dashV}{%
  \mathrel{%
    \mathpalette\@dashV\Vdash
  }%
}
\newcommand*{\@dashV}[2]{%
  \reflectbox{$\m@th#1#2$}%
}
\makeatother

\title{Topology as faithful communication\\ through relations}
\author{S.\,Maschio, G.\,Sambin}
\date{}

\begin{document}
\maketitle
\begin{abstract}
Basic pairs and their morphisms are the most elementary framework in which standard topological notions can be defined. We present here a new interpretation of topological concepts as those which can be communicated faithfully between the two sides of basic pairs. In particular, we prove that the subsets which can be communicated faithfully (in the suitable way) are exactly open subsets and closed subsets. We also prove that a relation (and in particular a function) between two sets of points can be communicated faithfully if and only if it is continuous.
\end{abstract}
\section{Introduction}
Most topological concepts can be presented in a predicative and constructive framework of basic pairs (see \cite{BP}). A basic pair $(X,\Vdash,S)$ consists of a set $X$, a set $S$ and a relation $\Vdash$ from $X$ to $S$. $X$ is called concrete side and represents points, while $S$ is called formal side and represents a set of indexes for a basis of neighbourhoods of a  topology on $X$. If $a$ is an index in $S$, then $\ext a$ is the subset of $X$ of those $x$ for which $x\Vdash a$. The presence of the formal side $S$ makes the structure of a basic pair \emph{symmetric}. If one adds two axioms
\begin{enumerate}
\item[$\mathsf{B1}$)] $\ext a\cap\ext b=\bigcup \{\ext c\,|\,\ext c\subseteq \ext a\cap \ext b\}$ 
\item[$\mathsf{B2}$)] $(\forall x\in X)(\exists a \in S)(x\Vdash a)$
\end{enumerate}
 one can obtain a predicative and constructive account of topological spaces.

Concepts in the concrete and formal side are obtained by moving information (subsets) from $X$ to $S$ or vice versa.

The idea is that topology can be read as faithful communication between its sides. 
This interpretation is not technically difficult, but it introduces a new intuitive point of view on topology which can shed light on unexpected links. 

In particular we show that the notions of open and closed sets can be interpreted as notions of communicable sets between the concrete and the formal side. In the same style we show that a relation is continuous if and only if it is communicable between the formal side and the concrete side. 
\subsection{Communication}
Suppose an individual $A$ wants to communicate with another individual $B$, but suppose $A$ and $B$ don't share the same language. 
However $A$ and $B$ both have their own collection of messages $M_{A}$ and $M_{B}$ which they use to represent information. Some messages in $M_{A}$ are equivalent, because they have the same meaning and the same holds for the messages in $M_{B}$. 
Such equivalences can be represented by equivalence relations $\sim_{A}$ and $\sim_{B}$ on $M_{A}$ and $M_{B}$ respectively.

Hence $A$ is equipped with a pair $(M_{A},\sim_{A})$ and $B$ with a pair $(M_{B},\sim_{B})$.

If we want $A$ and $B$ to communicate, then 
\begin{enumerate}
\item $B$ needs a decoding procedure $\Delta$ to transform every message in $M_{A}$ in one of its messages in $M_{B}$. This decoding procedure is good if it translates equivalent messages in $M_{A}$ in equivalent messages in $M_{B}$.
\item Conversely $A$ needs a decoding procedure $\nabla$ to transform every message in $M_{B}$ in one of its messages in $M_{A}$. This decoding procedure is good if it translates equivalent messages in $M_{B}$ in equivalent messages in $M_{A}$.
\end{enumerate}

As we well know, translators are not perfect devices and languages can be very different each others. 

We can say that a message $m$ in $M_{A}$ is (faithfully) communicable if it satisfies the following requirement: if $A$ communicates $m$ to $B$, $B$ translates it obtaining $\Delta(m)$  and then sends  $\Delta(m)$ back to $A$, then the translation $\nabla(\Delta(m))$ by $A$ of $\Delta(m)$ is equivalent to $m$.

We can make this idea precise in the following definition
\begin{definition}
A \emph{communication system} is a $4$-tuple $({\cal M}_{A},{\cal M}_{B},\Delta,\nabla)$ in which:
\begin{enumerate}
\item ${\cal M}_{A}$ is a pair $(M_{A},\sim_{A})$ with $M_{A}$ a collection and $\sim_{A}$ an equivalence relation on $M_{A}$
\item ${\cal M}_{B}$ is a pair $(M_{B},\sim_{B})$ with $M_{B}$ a collection and $\sim_{B}$ an equivalence relation on $M_{B}$
\item $\Delta$ is an operation from $M_{A}$ to $M_{B}$ (i.e.\ $\Delta(m)\in M_{B}\,[x\in M_{A}]$) such that for all $m,m'\in M_{A}$, if $m\sim_{A}m'$ then $\Delta(m)\sim_{B}\Delta(m')$
\item $\nabla$ is an operation from $M_{B}$ to $M_{A}$ (i.e.\ $\nabla(m)\in M_{A}\,[x\in M_{B}]$) such that for all $m,m'\in M_{B}$, if $m\sim_{B}m'$, then $\nabla(m)\sim_{A}\nabla(m')$.
\end{enumerate}
We say that $m\in M_{A}$ is $(\Delta,\nabla)$-communicable if $m\sim_{A}\nabla(\Delta(m))$. Similarly we say that $m'\in M_{B}$ is $(\Delta,\nabla)$-communicable if $m'\sim_{A}\Delta(\nabla(m'))$.
\end{definition}

\section{Communication of subsets}
\subsection{Operators on subsets in a basic pair}
We first recall some basic notions from \cite{BP}.

If $D$ and $E$ are subsets of a set $A$, then $D\between E$ is an abbreviation for the formula $(\exists x\in A)(x\releps  D\wedge x\releps  E)$.

If $\mathsf{r}$ is a relation from a set $X$ to a set $Y$, then if $x\in X$ and $y\in Y$ we can define a subset $\mathsf{r}\,x$ of $Y$ as $\{y\in Y|\,\mathsf{r}(x,y)\}$ and a subset $\mathsf{r}^{-}\,y$ of $X$ as $\{x\in X|\,\mathsf{r}(x,y)\}$.
If $D$ is a subset of $X$ and $E$ is a subset of $Y$ we define
\begin{equation}\label{r1}\mathsf{r}\,D:=\{y\in Y|\,\mathsf{r}^{-}\,y\between D\}\end{equation}
\begin{equation}\label{r2}\mathsf{r}^{-*}\,D:=\{y\in Y|\,\mathsf{r}^{-}\,x\subseteq D\}\end{equation}
\begin{equation}\label{r3}\mathsf{r}^{-}\,E:=\{x\in X|\,\mathsf{r}\,x\between E\}\end{equation}
\begin{equation}\label{r4}\mathsf{r}^{*}\,E:=\{x\in X|\,\mathsf{r}\,x\subseteq E\}\end{equation}

\begin{definition}
A \emph{basic pair} $(X,\Vdash, S)$ is a pair of sets $X$ and $S$ together with a relation $\Vdash$ from $X$ to $S$.
 \end{definition}
When we consider a basic pair we distinguish the operators defined in (\ref{r1}), (\ref{r2}), (\ref{r3}), (\ref{r4}) with a specific notation: 
\begin{definition}
Let $(X,\Vdash, S)$ be a basic pair. 
\begin{enumerate}
\item if $a\in S$, then $\ext a:=\,\Vdash^{-} a=\{x\in X|\,x\Vdash a\}$
\item if $x\in X$, then $\Diamond\,x:=\,\Vdash x=\{a\in S|\,x\Vdash a\}$
\end{enumerate}
If $D$ is a subset of $X$, then
\begin{enumerate}
\item $\Diamond D:=\,\Vdash D=\{a\in S|\,\ext a\between D) \}$
\item $\Box D:=\,\Vdash^{-*}D=\{a\in S|\,\ext a\subseteq D) \}$
\end{enumerate}
If $U$ is a subset of $S$ then
\begin{enumerate}
\item $\ext U:=\,\Vdash^{-} U=\{x\in X|\,\Diamond\,x\between U) \}$
\item $\rest \, U:=\,\Vdash^{*} U=\{x\in X|\,\Diamond\,x\subseteq U) \}$
\end{enumerate}
\end{definition}
Let us present some preliminaries results (see \cite{BP}).
\begin{prop}\label{prop}
The following properties hold for a basic pair $(X,\Vdash, S)$:
\begin{enumerate}
\item if $D\subseteq E$ are subsets of $X$, then $\Diamond D\subseteq \Diamond E$ and $\Box D\subseteq \Box E$;
\item if $U\subseteq V$ are subsets of $S$, then $\ext U\subseteq \ext V$ and $\rest \,U\subseteq \rest \,V$; 
\item $\ext$ is left adjoint to $\Box$, i.\,e.\ if $D$ is a subset of $X$ and $U$ is a subset of $S$, $\ext U\subseteq D$ if and only if $U\subseteq \Box D$;
\item $\rest $ is right adjoint to $\Diamond$, i.\,e.\ if $D$ is a subset of $X$ and $U$ is a subset of $S$, $D\subseteq\rest \,U$ if and only if $\Diamond D\subseteq U$.
\end{enumerate}
In particular if $D$ is a subset of $X$, then $\ext\Box D\subseteq D$ and $D\subseteq \rest \,\Diamond D$.
\end{prop}

\subsection{Communicable subsets}
Let $(X,\Vdash,S)$ be a basic pair. As we said in the introduction we call $X$ its concrete side and $S$ its formal side.

We are interested here in a communication system of the form

$$\left((\mathcal{P}X,=),(\mathcal{P}S,=),\Delta,\nabla\right)$$

The communicating individuals are the concrete side and the formal side of $(X,\Vdash,S)$. Their messages are subsets: ${\cal M}_{X}=(\mathcal{P}X,=)$ and ${\cal M}_{S}=(\mathcal{P}S,=)$ respectively.

The formal side can \emph{understand} a subset $D$ of $X$ of the concrete side only by means of a subset of $S$. There are basically two meaningful ways $\Delta$ for the formal side to understand the information in $D$: $\Box D$ or $\Diamond D$. The first is an approximation by defect of $D$: the formal side considers the subset of all neighbourhoods of which the extensions are contained in $D$. The second is an approximation by excess of $D$: the formal side considers the subset of all neighbourhoods of which the extensions overlap $D$. 

The same holds in the opposite direction. The concrete side can \emph{understand} a subset $U$ of $S$ of the formal side only by means of a subset of $X$. There are basically two meaningful ways $\nabla$ for the formal side to understand the information in $U$: $\rest \,U$ or $\ext U$. The first is an approximation by defect of $D$: the concrete side considers the subset of all points of which all basic neighbourhoods are in $D$. The second is an approximation by excess of $D$: the formal side considers the subset of all points which belong to some basic neighbourhood in $D$. 

We can interpret this communication in a different way. The formal side $S$ can \emph{understand} a subset $D$ of $X$ of the concrete side by asking the question ``Is $a$ in the concept $D$?'' to $X$. $X$ understands $a$ not by means of an element, but of the subset $\ext a$ and hence it can answer either ``Yes'' when $\ext a\between D$ or ``Yes'' when $\ext a\subseteq D$. In the first case $S$ understands $D$ as $\Diamond D$, in the second case $S$ understands $D$ as $\Box D$. A similar interpretation holds in the opposite direction.

\subsection{Open and closed subsets as communicables}
Let us fix a basic pair $(X,\Vdash,S)$.
\begin{definition}
A subset $D$ of $X$ is 
\begin{enumerate} 
\item \emph{open} if $(\forall x\in X)(x\releps  D\rightarrow (\exists a\in S)(x\Vdash a\wedge\ext a\subseteq D))$, i.e.\ if every point in $D$ has a basic neighbourhood included in $D$. 
\item \emph{closed} if $(\forall x\in X)((\forall a\in S)(x\Vdash a\rightarrow \ext a\between D)\rightarrow x\releps  D)$, i.e.\, if a point is in $D$ whenever every of its basic neighbourhoods overlaps $D$.
\end{enumerate}
\end{definition}

We first consider mixed decoding strategies ($(\Box, \ext)$ and $(\Diamond, \rest )$, which obviously make sense because if one side uses an approximation by defect, it is quite natural for the other to use an approximation by excess to compensate it.

Our first result is: 

\begin{theorem}\label{open}
Let $D$ be a subset of $X$. Then
\begin{enumerate}
\item $D$ is open if and only if $D$ is $(\Box,\ext)$-communicable;
\item $D$ is closed if and only if $D$ is $(\Diamond,\rest )$-communicable.
\end{enumerate}
\end{theorem}
\bpr By definition $D\textnormal{ is open}$ if and only if $(\forall x\in X)(x\releps  D\rightarrow \Diamond x\between \Box D)$ 
if and only if $D\subseteq\ext\Box D$. Similarly, by definition $D\textnormal{ is closed}$ if and only if $(\forall x\in X)(\Diamond x \subseteq \Diamond D\rightarrow x\releps  D)$ if and only if $\rest \,\Diamond D\subseteq D$. 

The result follows from proposition \ref{prop}.
$\Box$\\

In order to prove the next theorem we first need:
\begin{lemma}\label{extrest}
Let $U$ be a subset of $S$. Then $\ext U$ is open and $\rest U$ is closed.
\end{lemma}
\bpr
As $\ext U\subseteq \ext U$, from proposition it follows that \ref{prop} $U\subseteq \Box\ext U$ and hence $\ext U\subseteq \ext \Box\ext U$. We already know by proposition \ref{prop} that $\ext \Box\ext U\subseteq \ext U$. This implies that $\ext U=\ext \Box\ext U$. Hence, thanks to theorem \ref{open}, $\ext U$ is open.

Similarly one can prove that $\rest U$ is closed.
$\Box$\\

Let us now consider the excess/excess $(\Diamond, \ext)$ and defect/defect $(\Box, \rest )$ decoding strategies.

\begin{theorem}\label{de}
Let $D$ be a subset of $X$. Then 
\begin{enumerate}
\item if $D$ is $(\Diamond,\ext)$-communicable, then $D$ is open;
\item if $D$ is $(\Box ,\rest )$-communicable, then $D$ is closed.
\end{enumerate}
In particular if $D$ is $(\Box ,\rest )$-communicable and $(\Diamond ,\ext)$-communicable, then it is clopen.
\end{theorem}
\bpr
The proof is an immediate consequence of lemma \ref{extrest}.

$\Box$

%


The converses of the statements in theorem \ref{de} don't hold. 
In fact one can consider the basic pair $(\mathbf{2},\Vdash,\mathbf{3})$ where $x\Vdash y\equiv^{def} x=y\vee y=2$. 

The singleton subsets $\{0\}$ and $\{1\}$ of $\mathbf{2}$ are both closed and open, but none of them is either $(\Box ,\rest )$- or $(\Diamond ,\ext)$-communicable.

In particular the basic pairs $(\mathbf{2},\Vdash,\mathbf{3})$ and $(\mathbf{2},=,\mathbf{2})$ give rise to the same open and closed subsets of $\mathbf{2}$, while at the same time they give rise to different $(\Box ,\rest )$- and $(\Diamond ,\ext)$-communicable sets. \\


A basic pair $(X,\Vdash, S)$ satisfies the axiom $\mathsf{B2}$, if 
\begin{equation} X=\ext S \end{equation}

A useful consequence is 
\begin{lemma}\label{b2} If $(X,\Vdash, S)$ satisfies $\mathsf{B2}$, then for every subset $U$ of $S$
$$\rest \,U\subseteq \ext U$$
\end{lemma}
\bpr
Suppose $x\releps \rest \,U$. Then $\Diamond x\subseteq U$. As a consequence of $\mathsf{B2}$ there exists $a\in S$ such that $a\releps \Diamond x$. This implies that $a\releps U$.
We thus proved that $\Diamond x\between U$, i.\,e.\ $x\releps \ext U$.
$\Box$
\begin{theorem} If $(X,\Vdash, S)$ satisfies $\mathsf{B2}$, then for every subset $D$ of $X$ the following are equivalent:
\begin{enumerate}
\item $D$ is $(\Diamond,\ext)$-communicable;
\item $D$ is $(\Box,\rest )$-communicable;
\item $D$ is clopen and $\Diamond D\subseteq \Box D$.
\end{enumerate}
\end{theorem}
\bpr
\noindent $(1\Rightarrow 3)$ If $D$ is $(\Diamond,\ext)$-communicable, then $D=\ext\Diamond D$ and hence $\ext\Diamond D\subseteq D$. From this, by proposition \ref{prop}, it follows that $\Diamond D\subseteq \Box D$. Using proposition \ref{prop} and lemma \ref{b2} we can deduce the following chain of inclusions
$$ D\subseteq \rest \Diamond D \subseteq \ext\Diamond D\subseteq D$$
from which we can deduce that $D$ is closed. We already know that $D$ is open by theorem \ref{de}, thus $D$ is clopen.

\noindent $(2\Rightarrow 3)$
If $D$ is $(\Box,\rest)$-communicable, then $D=\rest\Box D$ and hence $D\subseteq \rest\Box D$. From this, by proposition \ref{prop}, it follows that $\Diamond D\subseteq \Box D$. Using proposition \ref{prop} and lemma \ref{b2} we can deduce the following chain of inclusions
$$ D\subseteq \rest \Box D \subseteq \ext\Box D\subseteq D$$
from which we can deduce that $D$ is open. We already know that $D$ is closed by theorem \ref{de}, thus $D$ is clopen.

\noindent $(3\Rightarrow 1),(3\Rightarrow 2)$ If $\Diamond D\subseteq \Box D$, then, using lemma \ref{b2}, we deduce that 
$$ \rest \Diamond D \subseteq \ext \Diamond D\subseteq \ext\Box D$$
and
$$ \rest \Diamond D \subseteq \rest \Box D\subseteq \ext\Box D$$
If $A$ is clopen, then $A=\ext\Box A=\rest\Diamond A$ and hence $A=\ext\Diamond A$ and $A=\rest\Box A$.
%

$\Box$

%
%

\subsection{Other decoding strategies}
There are other operators which one can define on the basic pair $(X,\Vdash,S)$:
\begin{definition}
If $D$ is a subset of $X$ and $U$ is a subset of $S$, then 
\begin{enumerate}
\item $D^{\rightarrow}:=\{a\in S|\,(\forall x\in X)(x\releps D\rightarrow x\Vdash a)\}=\{a\in S|\,D\subseteq \ext(a)\}$
\item $U^{\leftarrow}:=\{x\in X|\,(\forall a\in S)(a\releps U\rightarrow x\Vdash a)\}=\{x\in X|\,U\subseteq \Diamond U\}$
\end{enumerate}
\end{definition}
Notice that these operators enjoy the following properties:
\begin{prop}\label{lar}
If $D$ and $E$ are subsets of $X$ and $U$ and $V$ are subsets of $S$, then 
\begin{enumerate}
\item if $D\subseteq E$, then $E^{\rightarrow}\subseteq D^{\rightarrow}$;
\item if $U\subseteq V$, then $V^{\leftarrow}\subseteq U^{\leftarrow}$;
\item $D\subseteq U^{\leftarrow}$ if and only id $U\subseteq D^{\rightarrow}$;
\end{enumerate}
In particular $D\subseteq (D^{\rightarrow})^{\leftarrow}$.
\end{prop}
\begin{theorem}
Every $(^{\rightarrow},\ext)$-communicable subset of $X$ is open and every $(^{\rightarrow},\rest )$-communicable subset of $X$ is closed.
\end{theorem}
\bpr
This is an immediate consequence of lemma \ref{extrest}. $\Box$
\begin{theorem}
Every $(\Box,^{\leftarrow})$-communicable and every $(\Diamond,^{\leftarrow})$-communicable subset of $X$  is $(^{\rightarrow},^{\leftarrow})$-communicable.\end{theorem}
\bpr
If $D=(\Box D)^{\leftarrow}$, then in particular $D\subseteq(\Box D)^{\leftarrow}$ from which by proposition \ref{lar} it follows $\Box D\subseteq D^{\rightarrow}$ and thus $D\subseteq (D^{\rightarrow})^{\leftarrow}\subseteq(\Box D)^{\leftarrow}$.Hence $D=(D^{\rightarrow})^{\leftarrow}$.

Analogously one can prove that if $D=(\Diamond D)^{\leftarrow}$, then $D=(D^{\rightarrow})^{\leftarrow}$.
%
$\Box$
\begin{theorem}
$D$ is $(^{\rightarrow},^{\leftarrow})$-communicable if and only if $D=\bigcap\{\ext a|\,a\in U\}$ for some subset $U$ of $S$. 
\end{theorem}
\bpr
By definition 
$(D^{\rightarrow})^{\leftarrow}=
\bigcap\{\ext a|\,D\subseteq \ext a\}=\bigcap\{\ext a|\ a\releps D^{\rightarrow}\}$. 
Conversely, if $D=\bigcap\{\ext a|\,a\in U\}$, then $$(D^{\rightarrow})^{\leftarrow}=\bigcap\{\ext a|\, a\releps D^{\rightarrow}\}=\bigcap\{\ext a|\, D\subseteq \ext a\}=$$
$$=\bigcap\{\ext a|\,\bigcap\{\ext b|\,b\in U\}\subseteq \ext a\}\subseteq \bigcap\{\ext a|\,a\in U\}=D.$$ 
But we already know from lemma \ref{lar} that $D\subseteq (D^{\rightarrow})^{\leftarrow}$. Thus $D=(D^{\rightarrow})^{\leftarrow}$.
$\Box$
\begin{remark}
In a classical and impredicative framework a topology ${\cal T}$ on a set $\Omega$ gives rise to a basic pair $(\Omega, \in, {\cal T})$.
A $\mathsf{T}_{0}$-point is an equivalence class with respect to the equivalence relation of topological indistinguishibility $\sim$ on $\Omega$ defined by $x\sim y$ if and only if $\{D\in {\cal T}|\,x\in D\}=\{D\in {\cal T}|\,y\in D\}$. In this framework:
\begin{enumerate}
\item[]$D=\ext\Box D$ iff $D\in {\cal T}$.
\item[]$D=\rest \Diamond D$ iff $D$ is closed.
\item[]$D=\ext\Diamond D$ iff $D=\rest \Box D$ iff $D=\emptyset$ or $D=\Omega$.
\item[]$D=\ext(D^{\rightarrow})$ iff $D=((\Box D)^{\leftarrow})$ iff $D=\Omega$.
\item[]$D=((D^{\rightarrow})^{\leftarrow})$ iff $D$ is an intersection of open sets.
\item[]$D=\rest (D^{\rightarrow})$ iff $D$ is a closed $\mathsf{T}_{0}$-point.
\item[]$D=((\Diamond D)^{\leftarrow})$ iff $D$ is a $\mathsf{T}_{0}$-point which is an intersection of open sets.\end{enumerate}
\end{remark}
\section{Communication of relations}
\subsection{Relations between basic pairs}
Let us fix two basic pairs ${\cal X}=(X,\Vdash, S)$ and ${\cal Y}=(Y,\Vdash, T)$.

A relation $\mathsf{f}$ from $X$ to $Y$ is 
\begin{enumerate}
\item \emph{single-valued} if $(\forall x\in X)(\forall y\in Y)(\forall y'\in Y)(\mathsf{f}(x,y)\wedge \mathsf{f}(x,y')\rightarrow y=y')$;
\item \emph{total} if $(\forall x\in X)(\exists y\in Y)\,\mathsf{f}(x,y)$;
\item a \emph{function} if it is single-valued and total.
\end{enumerate}
As usual, if $\mathsf{f}$ is a function from $X$ to $Y$, for every proposition $P(y)$ depending on $y\in Y$ we define the abbreviation 
$$P(\mathsf{f}(x))\equiv^{def}(\exists y\in Y)(\mathsf{f}(x,y)\wedge P(y))$$
As $\mathsf{f}$ is a function, this is equivalent to $(\forall y\in Y)(\mathsf{f}(x,y)\rightarrow P(y))$.

Continuity of $\mathsf{f}$ with respect to ${\cal X}$ and ${\cal Y}$ can be defined, using basic neighbourhoods, as usual: for all $x\in X$ and for all $b\in T$
$$\mathsf{f}(x)\releps \ext b\rightarrow (\exists a\in S)(x\releps \ext a\wedge (\forall x'\in X)(x'\releps \ext a\rightarrow \mathsf{f}(x')\releps \ext b))$$
However one can easily notice that $P(\mathsf{f}(x))$ is equivalent to $\mathsf{f}x\between \{y\in Y|\,P(y)\}$.
Hence the condition of continuity is equivalent to: 
for all $x\in X$ and for all $b\in T$
$$\mathsf{f}x\between \ext b\rightarrow (\exists a\in S)(x\releps \ext a\wedge (\forall x'\in X)(x'\releps \ext a\rightarrow \mathsf{f}x'\between \ext b))$$
Using (\ref{r3}) the condition becomes: for all $x\in X$ and for all $b\in T$
$$x\releps \mathsf{f}^{-}\ext b\rightarrow (\exists a\in S)(x\releps \ext a\wedge (\forall x'\in X)(x'\releps \ext a\rightarrow x'\releps \mathsf{f}^{-}\ext b))$$
which is equivalent to: for all $x\in X$ and for all $b\in T$
$$x\releps \mathsf{f}^{-}\ext b\rightarrow (\exists a\in S)(x\Vdash a\wedge \ext a\subseteq \mathsf{f}^{-}\ext b)$$
We use this representation of continuity to extend the notion of continuity to relations.

\begin{definition}A relation $\mathsf{r}$ from $X$ to $Y$ is \emph{continuous} from ${\cal X}$ to ${\cal Y}$  if 
\begin{equation}\label{continuity}(\forall b\in T)(\forall x\in X)(x\releps  \mathsf{r}^{-}\ext b\rightarrow (\exists a\in S)(x\Vdash a\wedge \ext a\subseteq \mathsf{r}^{-}\ext b))\end{equation}
\end{definition}

Notice that the condition (\ref{continuity}) of continuity is also equivalent to:
\begin{equation}(\forall b\in T)(\forall x\in X)(x\releps  \mathsf{r}^{-}\ext b\rightarrow \Diamond x\between \Box\mathsf{r}^{-}\ext b))\end{equation}

We also define a notion of equivalence between relations from $X$ to $Y$ and from $S$ to $T$ with respect to ${\cal X}$ and ${\cal Y}$.
\begin{definition} For relations $\mathsf{r}_{1},\mathsf{r}_{2}$ from $X$ to $Y$, we write $\mathsf{r}_{1}\sim \mathsf{r}_{2}$ if 
$$(\forall b\in T)(\mathsf{r}_{1}^{-}\ext b=\mathsf{r}_{2}^{-}\ext b)$$
For relations $\mathsf{s}_{1},\mathsf{s}_{2}$ from $S$ to $T$, we write $\mathsf{s}_{1}\approx \mathsf{s}_{2}$ if 
$$(\forall x\in X)(\mathsf{s}_{1}\Diamond x=\mathsf{s}_{2}\Diamond\,x)$$
\end{definition}
\subsection{Communicable relations}
We are interested in a communication system of the form
$$\left((Rel(X,Y),\sim),(Rel(S,T),\approx),\Delta,\nabla\right)$$
Here the communicating individuals are the concrete sides $(X,Y)$ and the formal sides $(S,T)$. Their messages are relations: ${\cal M}_{(X,Y)}=(Rel(X,Y),\sim)$ and ${\cal M}_{(S,T)}=(Rel(S,T),\approx)$ respectively.

The next definition proposes candidates for $\Delta$ and $\nabla$. 

\begin{definition}
Let $\mathsf{r}$ be a relation from $X$ to $Y$, then $\sigma(\mathsf{r})$ from $S$ to $T$ is the relation defined by
$$\sigma(\mathsf{r})(a,b)\equiv^{def}\ext a\subseteq \mathsf{r}^{-}\ext b$$
Let $\mathsf{s}$ be a relation from $S$ to $T$, then $\rho(\mathsf{s})$ from $X$ to $Y$ is the relation defined by
$$\rho(\mathsf{s})(x,y)\equiv^{def}\Diamond y\subseteq \mathsf{s}\Diamond x$$

\end{definition}

Let us first prove the following:
\begin{prop}
Let $\mathsf{r}_{1},\mathsf{r}_{2}$ be relations from $X$ to $Y$ and let $\mathsf{s}_{1},\mathsf{s}_{2}$ be relations from $S$ to $T$. Then
\begin{enumerate}
\item if $\mathsf{r}_{1}\sim\mathsf{r}_{2}$, then $\sigma(\mathsf{r}_{1})=\sigma(\mathsf{r}_{2})$;
\item if $\mathsf{s}_{1}\approx\mathsf{s}_{2}$, then $\rho(\mathsf{s}_{1})=\rho(\mathsf{s}_{2})$.
\end{enumerate} 
\end{prop}
\bpr
if $\mathsf{r}_{1}\sim\mathsf{r}_{2}$, then for every $b\in T$, $ \mathsf{r}_{1}^{-}\ext b=\mathsf{r}_{2}^{-}\ext b$ and thus
$$\sigma(\mathsf{r}_{1})(a,b)\leftrightarrow \ext a\subseteq \mathsf{r}_{1}^{-}\ext b\leftrightarrow \ext a\subseteq \mathsf{r}_{2}^{-}\ext b\leftrightarrow\sigma(\mathsf{r}_{2})(a,b)$$

The proof of the second statement is analogous.
$\Box$\\

Hence we obtain the following:
\begin{cor}
$((Rel(X,Y),\sim),(Rel(S,T),\approx),\rho,\sigma)$ is a communication system.
\end{cor}


\subsection{Continuity as communication}
Before proving the main theorem of this section, we need some preliminary lemmas:
\begin{lemma} Let $\mathsf{r}$ be a relation from $X$ to $Y$. Then for all $b\in T$
$$\rho(\sigma(\mathsf{r}))^{-}\ext b\subseteq \mathsf{r}^{-}\ext b$$
\end{lemma}
\bpr
Let $b\in T$ and suppose $x\releps \rho(\sigma(\mathsf{r}))^{-}\ext b$.
Thus there exists $y\in Y$ such that $y\Vdash b$ and $\rho(\sigma(\mathsf{r}))(x,y)$. Notice that $\rho(\sigma(\mathsf{r}))(x,y)$ is equivalent to $$(\forall c\in T)(c\releps \Diamond y\rightarrow \Diamond x\between \Box \mathsf{r}^{-}\ext c)$$

As $b\releps \Diamond y$ we obtain $\Diamond x\between \Box \mathsf{r}^{-}\ext b$ which is equivalent by definition to $x\releps \ext \Box \mathsf{r}^{-}\ext b$ from which it immediately follows by proposition \ref{prop} that $x\releps \mathsf{r}^{-}\ext b$.

\begin{lemma} If $\mathsf{r}:X\rightarrow Y$ is continuous, then for all $b\in T$
$$\mathsf{r}^{-}\ext b\subseteq \rho(\sigma(\mathsf{r}))^{-}\ext b$$
\end{lemma}
\bpr
Suppose $\mathsf{r}$ is continuous. Then for every $x\in X$ and for every $c\in T$
$$x\releps \mathsf{r}^{-}\ext c\rightarrow \Diamond x\between \Box \mathsf{r}^{-}\ext c$$
Moreover suppose that $x\releps \mathsf{r}^{-}\ext b$. This implies that we can fix an $y\in Y$ such that $b\releps \Diamond y$ and $\mathsf{r}(x,y)$. 
Suppose now that $c\releps \Diamond y$ and $c\in T$. This implies that $y\releps \mathsf{r}^{-}\ext c$ and hence, by continuity,  $\Diamond x\between \Box \mathsf{r}^{-}\ext c$.
So we proved that 
$$(\exists y\in Y)(b\releps \Diamond y\wedge (\forall c\in T)(c\releps \Diamond y\rightarrow \Diamond x\between \Box \mathsf{r}^{-}\ext c ))$$
i.e.\ $x\,\releps \rho(\sigma(\mathsf{r}))^{-}\ext b$. $\Box$

%
%
%
%

\begin{theorem}
$\mathsf{r}$ is continuous if and only if $\mathsf{r}$ is $(\sigma,\rho)$-communicable.
\end{theorem}
\bpr 
If $\mathsf{r}$ is continuous, then $\mathsf{r}\sim\rho(\sigma(\mathsf{r}))$ as it follows from the previous two lemmas.

Conversely if  $\mathsf{r}\sim \rho(\sigma(\mathsf{r}))$, then in particular for every $x\in X$ and $b\in T$
$$(x\releps\mathsf{r}^{-}\ext b\rightarrow x\releps \rho(\sigma(\mathsf{r}))^{-}\ext b)$$

Suppose now that $x\releps\mathsf{r}^{-}\ext b$, then $x\releps \rho(\sigma(\mathsf{r}))^{-}\ext b$ which in particular means that there exists $y\in Y$ with $y\Vdash b$ and $ \rho(\sigma(\mathsf{r}))(x,y)$, i.e.\,
$$(\forall c\in T)(c\releps\Diamond y\rightarrow \Diamond x\between \Box \mathsf{r}^{-}\ext c)$$
Taking $b=c$ and using the fact that $b\releps \Diamond y$ one obtains that  $\Diamond x\between \Box \mathsf{r}^{-}\ext b$.

Thus $x\releps\mathsf{r}^{-}\ext b\rightarrow \Diamond x\between \Box \mathsf{r}^{-}\ext b$. Hence $\mathsf{r}$ is continuous. $\Box$\\

Let us now consider the case in which $\mathsf{r}$ is a function.

A basic pair $(X,\Vdash, S)$ is \emph{Hausdorff} (or $\mathsf{T}_{2}$) if for every $x\in X$ and $x'\in X$
$$(\forall a\in S)(\forall a'\in S)(x\Vdash a \wedge x'\Vdash a'\rightarrow \ext a\between \ext a')\rightarrow x=x'$$

\begin{prop}
If $\mathsf{f}$ is a function and $(Y,\Vdash, T)$ is Hausdorff, then $\rho(\sigma(\mathsf{f}))$ is a single-valued relation which is a restriction of $\mathsf{f}$.
\end{prop}
\bpr
Suppose that $\rho(\sigma(\mathsf{f}))(x,y)$. If $y\Vdash b$, then there exists $a\in S$, such that $x\Vdash a$ and $\ext a\subseteq \mathsf{f}^{-}\ext b$, and thus in particular 
 $\mathsf{f}(x)\Vdash b$.
This implies that if $y\Vdash b$ and $\mathsf{f}(x)\Vdash b'$, then $\ext b\between \ext b'$. As $(Y,\Vdash, T)$ is Hausdorff, we obtain that $y=\mathsf{f}(x)$.
Hence $\rho(\sigma(\mathsf{f}))$ is a restriction of $\mathsf{f}$ and thus it is a single-valued relation. $\Box$
\bibliographystyle{plain}
\bibliography{bibliopsp}
\end{document}